\documentclass{ems-mag-author}

\usepackage{blindtext}
\usepackage{tikz}
\usepackage{graphicx}
\usetikzlibrary{patterns}
\theoremstyle{plain}
\newtheorem{theorem}{Theorem}[section]

\theoremstyle{definition}

\newtheorem*{remark}{Remark}

\newlength{\defbaselineskip} 
\newcommand{\mint}{\mathop{\int\hskip -1,05em -\, \!\!\!}\nolimits}

\begin{document}

\doi{???}
\issue{???}

\def\eqn#1$$#2$${\begin{equation}\label#1#2\end{equation}}
\def\dx{\,{\textnormal{d}}x}
\def\dy{\,{\textnormal{d}}y}
\newcommand{\snr}[1]{\lvert #1\rvert}
\newcommand{\loc}{\textnormal{loc}}
\newcommand{\Dim}{\textnormal{dim}_{\mathcal{H}}}
\newcommand{\nr}[1]{\lVert #1 \rVert}
\newcommand{\tx}[1]{\textnormal{\texttt{#1}}}
\def\BBB{\textnormal{B}}
\title{\textbf{$\mu$}-ellipticity and nonautonomous integrals}

\author{Cristiana De Filippis}

\maketitle

\begin{abstract}
$\mu$-ellipticity is a form of nonuniform ellipticity arising in various contexts from the calculus of variations. Understanding regularity properties of minimizers in the nonautonomous setting is a challenging task fostering the development of delicate techniques and the discovery of new irregularity phenomena. 
\end{abstract}

\noindent
The classical area functional, given by
\begin{equation}\label{0.3}
    w\mapsto \int_{\Omega}\sqrt{1+\snr{Dw}^{2}}\dx
\end{equation}
and its Euler-Lagrange equation, the celebrated Minimal Surface Equation
\begin{equation}\label{0.0}
    -\textnormal{div}\left(\frac{Du}{\sqrt{1+\snr{Du}^2}}\right)=0\qquad \mbox{in} \ \ \Omega
\end{equation}
are classical objects of study in the modern calculus of variations and in theory of elliptic partial differential equations.\footnote{Unless otherwise specified in this note we shall always assume that $\Omega\subset \mathbb R^n$ is a bounded open set and $n\geq 2$. Moreover, we shall denote $B_{\varrho}(x)=\{y\in\mathbb R^n \, \colon \,\snr{y-x}<\varrho \}$.} Their peculiarities allowed to build a rich and large existence and regularity theory and have fostered generations of mathematicians to tackle difficult analytical questions. Equation \eqref{0.0} is intimately linked to the classical Plateau problem, which has historically driven the development of Geometric Measure Theory. Foundational contributions by De Giorgi, Reifenberg, Federer, Fleming, Almgren, Simons, Bombieri, Miranda, and Giusti have shaped the field. In this note we are particularly interested in gradient estimates for solutions and minimizers of integral functionals featuring ellipticity properties connected to the ones of \eqref{0.3}. As for \eqref{0.3}, we'd like to single out a particularly elegant result of Bombieri, De Giorgi and Miranda \cite{bdm69}, see also Trudinger \cite{tr72}, asserting the validity of the pointwise gradient estimate 
\begin{equation}\label{0.2}
    \snr{Du(x)}\lesssim_n\exp\left(c(n)\sup_{y\in B_{\varrho}(x)}\frac{\snr{u(y)-u(x)}}{\varrho}\right),
\end{equation}
for $C^2$-solutions $u$ to \eqref{0.0}. This a priori estimate is an essential tool in the proof of existence theorems of classical solutions, see \cite[Theorem 13.6]{giu84}. The functional in \eqref{0.3} is an example of a general integral of the calculus of variations of the type 
\begin{equation}\label{fun}
w   \mapsto \mathcal{F}(w,\Omega):= \int_{\Omega} \tx{F}(x, Dw)\dx
\end{equation}
where $\tx{F}\colon \Omega \times \mathbb R^n \to \mathbb R$ is a Carath\'eodory integrand\footnote{That is $x \mapsto \tx{F}(x,z)$ is measurable for every $z$ and $z \mapsto \tx{F}(x,z)$ is continuous for a.e. $x$. This ensures that the composition $x \mapsto \tx{F}(x, \tx{D}(x))$ is measurable whenever $\tx{D}\colon\Omega \to  \mathbb R^n$  is a measurable vector field.} having linear growth, in the sense that $\tx{F}(x, Dw)\approx \snr{Dw}$ for $\snr{Dw}$ large, see \cite{bs13}. Another such example is provided by the integral 
\begin{equation}\label{0.6.2}
w\mapsto \int_{\Omega}\left(1+\snr{Dw}^{m}\right)^{1/m}\dx,\qquad m>1\,.
\end{equation}
Next, consider the superlinear, $p$-growth classical model 
\begin{equation}\label{plap}
w\mapsto \int_{\Omega}(1+\snr{Dw}^{2})^{p/2}\dx, \qquad    p>1\,.
\end{equation}
Also in this case we have a neat a priori gradient estimate for minimizers $u$, i.e., 
\begin{equation*}
    \snr{Du(x)}\lesssim_{n,p} \Big( \frac{1}{\snr{B_{\varrho}(x)}}\int_{B_{\varrho}(x)} \snr{Du}^p\dy \Big)^{1/p}+1
\end{equation*}
which can be derived as in the fundamental work of Ural'tseva \cite{ura68} and Uhlenbeck \cite{uhl77}. Functionals with superlinear growth as in \eqref{plap} are at the core of a vast part of by now classical literature. Here we shall mainly concentrate on a class of  borderline integrals lying in between those with linear growth as in \eqref{0.3} and \eqref{0.6.2} and those with standard polynomial growth as in \eqref{plap}. These are functionals as \eqref{fun} with so-called nearly linear growth, i.e., such that 
 \begin{equation}\label{crescite}
 \displaystyle 
 \lim_{\snr{z}\to \infty} \frac{\tx{F}(x,z)}{\snr{z}}=\infty , \quad 
\lim_{\snr{z}\to \infty} \frac{\tx{F}(x,z)}{\snr{z}^{p}}=0 \quad \mbox{for all} \ \ p>1.
\end{equation}
A typical example belonging to such a class is the $L\log L$ functional 
\begin{equation}\label{0.6}
w\mapsto \int_{\Omega}\snr{Dw}\log(1+\snr{Dw})\dx,
\vspace{0.2mm}
\end{equation}
and its iterated versions 
\begin{equation}\label{0.7}
w\mapsto \int_{\Omega}\tx{L}_{i+1}(Dw)\dx,
\vspace{0.2mm}
\end{equation}
where, for integer $i\ge 0$, the integrands $\tx{L}_{i+1}$ are inductively defined via
\begin{equation}\label{0.6.1}
         \begin{cases}
         \displaystyle
    \ \ell_{0}(\snr{z})=\snr{z}\vspace{1mm}\\
    \displaystyle
    \ \ell_{i+1}(\snr{z}):=\log(1+\ell_{i}(\snr{z}))\quad &\mbox{for} \ \ i\ge 0\vspace{1mm}\\
    \displaystyle
    \ \texttt{L}_{i+1}(z):=\snr{z}\ell_{i+1}(\snr{z})\quad &\mbox{for} \ \ i\ge 0,
    \end{cases}
\end{equation}
for all $z\in \mathbb{R}^{n}$. As a consequence of the superlinear growth in \eqref{crescite}, the functionals we shall consider in the following pages will always be defined on the Sobolev space $W^{1,1}$, where in this situation Direct Methods of the calculus of variations apply. Indeed, we shall consider situations where 
\begin{equation}\label{superl}
\tx{G}(\snr{z})\lesssim \tx{F}(x,z) \,,\quad  \mbox{where} \lim_{\snr{z}\to \infty}\frac{\tx{G}(\snr{z})}{\snr{z}}= \infty\,,
\end{equation} which implies the first condition in \eqref{crescite}; this allows to recover weak compactness in $W^{1,1}$ of minimizing sequences via classical Dunford-Pettis's theorem. This is for instance the case of \eqref{plap}, \eqref{0.6} and \eqref{0.7}. Accordingly, a function $u \in W^{1,1}_{\loc}(\Omega)$ will be called a local minimizer\footnote{From now on, simply a minimizer.} of the functional $\mathcal F$ if for every ball $\BBB\Subset \Omega$ we have $
\tx{F}(\cdot, Du) \in L^1({\BBB}) $ and 
 $\mathcal{F}(u,\BBB)\leq \mathcal{F}(w,{\BBB})$  holds for every $w \in u + W^{1,1}_0({\BBB})$.

\section{Anisotropic $\mu$-ellipticity}\label{mumu1}
In view of \eqref{crescite}, a natural way to quantify the ellipticity properties of the functional in \eqref{fun}, and in such a way to cover all the models considered above, is to use the concept of (anisotropic) $\mu$-ellipticity. We assume that  $z\mapsto \tx{F}(\cdot, z)$ is $C^2$-regular and satisfies \vspace{0.2mm}
\begin{equation}\label{0.5}
\ \ \ \ \ \ \frac{\snr{\xi}^{2}}{(\snr{z}^{2}+1)^{\mu/2}} \lesssim \langle\partial_{zz}\tx{F}(x,z)\xi,\xi\rangle\lesssim \frac{(1+\texttt{g}(\snr{z}))\snr{\xi}^{2}}{(\snr{z}^{2}+1)^{(2-q)/2}}
\end{equation}
for all $z,\xi\in \mathbb{R}^{n}$, $x\in \Omega$, where $\mu\in (2-q,\infty)$, $q\ge 1$ are fixed numbers and $\tx{g}\colon [0,\infty)\to (0,\infty)$ is a continuous, nondecreasing, possibly unbounded function, with at most power growth at infinity. Related equations of the type $
-\textnormal{div} \, A(x,Du)=0
$
arising in connection with the Euler-Lagrange equation of the functional in \eqref{fun}, i.e., 
\begin{equation}\label{eulero}
-\textnormal{div} \, \partial_z \tx{F}(x,Du)=0,
\end{equation}
can be also considered. In this case we can use assumptions \eqref{0.5}  with $\partial_{zz}\tx{F}$ replaced by $\partial_{z} A$, where $A\colon \mathbb R^n\mapsto \mathbb R^n$ is a $C^1$-vector field.
Note that the integrand appearing in \eqref{0.3} fits \eqref{0.5} with $\mu=3$, $q\geq 1$, $\tx{g}(\snr{z})\equiv 1$, while the one in \eqref{0.6} verifies \eqref{0.5} with $\mu=1$, $q \geq 1$ and $\tx{g}(\snr{z})\equiv \log(1+\snr{z})$\footnote{The integrand in \eqref{0.6.2}  verifies (for $\snr{z}\not= 0$ when $m<2$)
\begin{equation*}
\scalebox{0.98}{$
\ \  \frac{\min\{m-1,1\}|z|^{m-2} \snr{\xi}^{2}}{(1+|z|^m)^{2-1/m}} \leq \langle\partial_{zz}\tx{F}(z)\xi,\xi\rangle\leq \frac{\max\{m-1,1\}|z|^{m-2} \snr{\xi}^{2}}{(1+|z|^m)^{1-1/m}}$}
\end{equation*} 
so that  \eqref{0.5} are satisfied with $\mu=m+1$, $q=1$ and $\tx{g}(\snr{z})\equiv 1$ provided $\snr{z}\geq \delta>0$ (the constants involved in \eqref{0.5} depend on $\delta$). Functionals as in \eqref{0.6.2} are studied for instance \cite{bs15, parks}.
}.The integrands $\tx{L}_{i+1}$ in \eqref{0.7}-\eqref{0.6.1} instead satisfy \eqref{0.5} for  any $\mu>1$, $q=1$ and with $\tx{g}(\snr{z})\equiv \ell_{i+1}(\snr{z})$, cf. \cite{fs99,fm00,dm23b,ddp24}. 
 Finally, \eqref{plap} satisfies \eqref{0.5} with $\mu=2-p$, $q=p>1$ and $g(\snr{z})\equiv  1$. Functionals \eqref{0.6}-\eqref{0.7} appear for instance in the theories of Prandtl-Eyring fluids and plastic materials with logarithmic hardening, \cite{fs99}, see also \cite{bil03} for more examples and a detailed discussion. The Orlicz space $L\log L(\Omega)$, defined via $f \in L\log L (\Omega)$ if and only if $\snr{f}\log(1+\snr{f})\in L^1(\Omega)$, directly connects to the functional in \eqref{0.6} and plays a crucial role in modern analysis, especially for its relations to Hardy spaces and maximal operators \cite{stein}. 
 
 \section{Nonuniform ellipticity and degeneracy} $\mu$-ellipticity is a degenerate type of nonuniform ellipticity in the sense that the lowest eigenvalue of $\partial_{zz}\tx{F}$ might, in principle, admit no positive lower bound. This follows considering the so-called ellipticity ratio, defined as
\begin{equation}\label{0.8}
\mathcal{R}_{ \tx{F}}(x,z):=\frac{\mbox{highest eigenvalue of} \ \partial_{zz}\tx{F}(x,z)}{\mbox{lowest eigenvalue of} \ \partial_{zz}\tx{F}(x,z)}\,.
\vspace{0.2mm}
\end{equation}
The boundedness of such a quantity is the condition defining classical  uniform ellipticity for  equations and functionals \cite{simon}, and, in that setting, it is crucial to derive a priori estimates for solutions. Here the situation is different. Condition \eqref{0.5} yields that the only a priori available bound on the ellipticity ratio is
\begin{equation}\label{0.9}
\mathcal{R}_{\tx{F}}(x,z)\lesssim \tx{g}(\snr{z})\snr{z}^{\mu+q-2}, \quad \mbox{for $\snr{z}\ge 1$}
\end{equation}
that yields no uniform control for $\snr{z}\to \infty$ (when $\mu+q>2$ and when $\mu+q=2$ and $g(\snr{z})\to\infty$ for $\snr{z}\to \infty$). This occurrence pushes $\mu$-elliptic problems out of reach for regularity techniques of standard use in the uniformly elliptic setting \cite{ura68,uhl77,gg82,gg83,gg84}. Degeneracy represents another pathological feature of $\mu$-ellipticity. As indicated by \eqref{0.5}, the smallest eigenvalue of $\partial_{zz} \tx{F}$ - which characterizes the ellipticity of the operator - has a power-type decay at infinity with respect to the gradient variable $z$. This, due to severe loss of ellipticity, makes the regularity theory of $\mu$-elliptic problems very challenging, rich and technically delicate. 

\section{The autonomous case}\label{mumu}
 The first general gradient regularity result for general $\mu$-elliptic integrals available in the literature is the following theorem.
\begin{theorem}[Fuchs and Mingione \cite{fm00}]\label{t2}
 Let $u\in W^{1,1}_{\textnormal{loc}}(\Omega)$ be a minimizer of functional \eqref{fun} with $\tx{F}(x,z)\equiv \tx{F}(z)$ 
satisfying $\tx{F}(z)\lesssim \snr{z}^q+1$, \eqref{superl} and \eqref{0.5} with $g(\cdot)\equiv 1$, $q \in (1,2)$, $1\le \mu<2$ such that
 \begin{equation}\label{0.12}
\mu+q<2+2/n.
 \end{equation}
Then $u\in W^{1,\infty}_{\textnormal{loc}}(\Omega)$. 
\end{theorem}
\noindent This covers the models $\tx{L}_i$ in \eqref{0.6.1}. In most of the results on nonuniformly elliptic problems we shall consider, the core point is actually to prove that $Du$ is locally bounded. Once this is secured, we are in a sense back to the uniformly elliptic setting,\footnote{Indeed, by \eqref{0.9} the ellipticity ratio $\mathcal{R}_{\tx{F}}(x,Du)$ cannot blow up when $Du\in L^\infty$. Therefore the problem behaves as it was uniformly elliptic when considered on Lipschitz solutions.} and more standard, yet delicate methods can be adapted to obtain higher regularity of minima, and, in particular, local H\"older continuity of first derivatives of minima (see \cite[Section 10]{dm23a},  \cite[Sections 5.9-5.11]{dm23b}, \cite[Section 5]{dm25}). In the case of Theorem \ref{t2}, the local H\"older continuity of $Du$ follows, i.e.,
\eqn{classicalsc0}
$$
\mbox{$C^{0,1}$-estimates $\Longrightarrow C^{1,\beta}$-estimates}
$$
and, in the case of Theorem \ref{t2}, for every $\beta \in (0,1)$. 
In view of \eqref{0.9}, conditions of the type in 
\eqref{0.12} obviously limit the growth of the ellipticity ratio $\mathcal{R}_{\tx{F}}(\cdot,Du)$ with respect to $Du$, while in fact proving that the gradient is locally bounded. In order to enlarge the rate of nonuniform ellipticity of the problem considered, that is, to allow a larger value of  $\mu+q$, it is possible to incorporate interpolative information, such as, for instance, some a priori boundedness on solutions.\footnote{This is for instance implied by maximum principle, when minimizers are found solving Dirichlet problems with bounded boundary data.} This has the ultimate effect of dropping the dimensional dependence on the growth of the ellipticity ratio, i.e.,  \eqref{0.12} can be replaced by 
    \begin{equation}\label{0.13}
    \mu +q<4,  \qquad u \in L^{\infty}_{\loc}(\Omega).
    \end{equation} 
\noindent For results of this type we refer to \cite[Section 5.2]{bil03}. Theorem \ref{t2} rests on an anisotropic version of Moser-type iteration, whose convergence is ensured by \eqref{0.12}. In the \eqref{0.13}-variant case, this method also involves a careful use of certain interpolation-type inequalities aimed at maximizing the integrability gain, eventually leading to the relaxed bound in \eqref{0.13}.  In both cases, the first step of the proof consists in differentiating the Euler-Lagrange equation \eqref{eulero},
which unavoidably breaks down when considering nonautonomous integrands with nondifferentiable coefficients, like for instance, when $x\mapsto \partial_z\tx{F}(x, \cdot)$ is only H\"older continuous. In this case, scheme \eqref{classicalsc0} is not viable using standard methods and novel ideas  must be developed, as we shall see in the next sections. For further literature on the autonomous case we recommend the interesting work of Marcellini and Papi \cite{marpapi}.

\begin{remark}[Vectorial problems] 
In this note we deal with the scalar case, i.e., when minima and competitors are scalar functions. Nevertheless, a large literature is available on the vectorial one, depending on the kind of regularity one is interested in. In general, and already in the uniformly elliptic case, solutions to elliptic systems and minima of vectorial functionals might exhibit singularities even in the most favourable situation of smooth, autonomous integrands. What is usually done in those cases is proving partial regularity, i.e., regularity of minima outside a negligible closed subset whose Hausdorff dimension can be eventually proven to be smaller than the ambient dimension; we refer to \cite{giabook} for an account of this theory. Additional structural assumptions on the integrand allow to prove everywhere regularity in the interior. For instance, Theorem \ref{mumu} extends to vector valued solutions 
 provided $\tx{L}$ is assumed to have a so-called Uhlenbeck structure \cite{uhl77}, i.e.,  $\tx{L}(Dw)=\ell (\snr{Dw})$. Partial regularity results in the vectorial case under $\mu$-ellipticity conditions  were established by Bildhauer and Fuchs \cite{bf01}. Key advancements are due to Gmeineder and Kristensen \cite{gk24}, who developed a unified, sharp approach to the almost everywhere regularity of minima of anisotropic multiple integrals covering also nonconvex, possibly signed functionals; see also \cite{gme21,gk19, sch09} for earlier results, and \cite{def22} where optimal partial regularity criteria are inferred via Nonlinear Potential Theory. Finally, dimensionless bounds as $q<p+2$ have been employed in the vectorial case in \cite{ckp11} by means of certain tricky penalization methods. 
\end{remark}

\subsection{Superlinear nonuniform ellipticity}\label{pqpq} \vspace{1.5mm}The bounds relating the size of $\mu$ and  $q$ in \eqref{0.12} and \eqref{0.13} are the natural counterpart of those appearing in the theory of nonuniformly elliptic problems with superlinear growth. Specifically, they parallel those available for so-called functionals with $(p,q)$-nonuniform ellipticity \cite{mar89,mar91}, formulated as 
\begin{equation}\label{0.14}
\frac{\snr{\xi}^{2}}{(\snr{z}^{2}+1)^{(2-p)/2}}\lesssim \langle\partial_{zz}\tx{F}(x,z)\xi,\xi\rangle\lesssim \frac{\snr{\xi}^{2}}{(\snr{z}^{2}+1)^{(2-q)/2}}
\end{equation}
for all $x\in \Omega$, $z,\xi\in \mathbb{R}^{n}$ and exponents $1<p\leq q$. Accordingly, $(p,q)$-growth conditions refer to similar conditions but this time prescribed directly on the integrand, i.e., 
\eqn{soddi}
$$
\snr{z}^p - c \lesssim \tx{F}(x, z) \lesssim \snr{z}^q+1\,, \quad c \geq 0\,.
$$
Conditions \eqref{0.14} and \eqref{soddi} are often verified together when considering autonomous, convex integrands. 
In such situations uniform ellipticity in ensured only when $p=q$. Formally, conditions \eqref{0.5} and \eqref{0.14} coincide letting $\mu=2-p$. Restrictions on the size of the so-called gap $q/p$ are necessary for minima to be regular. 
\begin{theorem}[Giaquinta \cite{gia87}; Marcellini \cite{mar87,mar89}] Let $\Omega\subset \{x\in \mathbb R^n \colon x_n>0\}$ be an open, bounded set. With $n>3$ and $q>2$, the function\footnote{We denote points  $x\in \mathbb R^n$ as $x=(x_1, \dots, x_n)$.}
\begin{equation}\label{minima}
u(x):=\left(\frac{c_{n,q}x_{n}^{q}}{\sum_{i=1}^{n-1}\snr{x_i}^{2}}\right)^{\frac{1}{q-2}},
\end{equation}
where
\begin{equation*}
c_{n,q}:=\left(\frac{n-1}{q-1}-\frac{2}{q-2}\right)\left(\frac{q-2}{q}\right)^{q-1}, 
\end{equation*}
is a minimizer of the functional 
\begin{equation}\label{0.15}
w\mapsto \int_{\Omega}\left(\frac{1}{2}\sum_{i=1}^{n-1}\snr{D_{i}w}^{2}+\frac{1}{q}\snr{D_{n}w}^{q}\right)\dx
\end{equation}
provided
\begin{equation}\label{0.16}
q>\frac{2(n-1)}{n-3}.
\end{equation}
\end{theorem}
\noindent The integrand in \eqref{0.15} satisfies condition \eqref{soddi} with $p=2$ and the function $u$ is obviously unbounded on the line $(0,\cdots,0,x_{n})$.
Similar examples can be produced with functionals satisfying \eqref{0.14} with $p=2$, see \cite{hong} and \cite{mar89}, thus offering instances of convex, scalar, regular integrals, with nonsmooth minimizers. This stands in sharp contrast with the classical literature, when in the case $p=q$ solutions and minimizers typically have H\"older continuous gradient \cite{ura68}.  On the positive side, in violation of \eqref{0.16} with $p=2$, we have the following theorem.
\begin{theorem}[Hirsch and Schäffner \cite{hs21}]\label{hst}
Let  $u \in W^{1,1}_{\loc}(\Omega)$ be a minimizer of the functional
\eqref{fun}, 
where the autonomous integrand $\tx{F}\colon \mathbb R^n \to \mathbb R$ is strictly convex and satisfies  \eqref{soddi} with 
\begin{equation}\label{ottimale}1 < p \leq q , \qquad \frac 1p-\frac 1q \leq  \frac{1}{n-1}\,.
\end{equation}
Then $u \in L^{\infty}_{\loc}(\Omega)$\footnote{In \cite{hs21}, Hirsch and Schäffner consider nonautonomous  Carath\'eodory integrands $\tx{F}\colon \Omega \times \mathbb R^n \to \mathbb R$, and assume convexity of $z \mapsto \tx{F}(\cdot, z)$ and  $\tx{F}(\cdot, 2z)\lesssim \tx{F}(\cdot, z)$. This is in line with the traditional De Giorgi-Nash-Moser theory, where, for the level of regularity of solutions considered here, coefficients can be allowed to be just measurable. The version reported here can be obtained by combining the a priori estimates of Hirsch and Schäffner with an approximation argument as for instance the one in \cite[Section 8]{dm25}.}. 
\end{theorem}
\noindent Theorem \ref{hst} builds on earlier results of Bella \& Sch\"affner \cite{bs20}. Note that in this case no assumptions on second derivatives of the integrand $\tx{F}$ of the type in \eqref{0.5} are imposed and nonuniform ellipticity is implicitly described by the growth conditions in \eqref{soddi}. Accordingly, no gradient regularity of minima is involved. The interest in the previous result, making it standing out in comparison to the previously published literature, rests on the optimal condition on the exponents $(p,q)$ in \eqref{ottimale}. Similar general sharp results remain unknown when switching to gradient regularity and considering assumptions \eqref{0.14}, while positive results on gradient boundedness are in \cite{bs20, mar91}. However, for certain large classes of functionals it is possible to derive sharp bounds on $q/p$, as shown in the work of  Koch, Kristensen and the author \cite{dkk24}, who cover autonomous integrands $\tx{F}(x,Dw)\equiv \tx{F}(Dw)$ that are convex, even polynomials, with non-negative homogeneous components and lowest homogeneity degree larger than $p\geq 2$. Indeed, the peculiar structure of convex polynomials allows for a finer nonuniform ellipticity measurement, referred to in \cite{dkk24} as Legendre $(p,q)$-nonuniform ellipticity, which quantifies the subtle interplay between the gradient of minima and the stress tensor. This is analysed  via convex duality arguments and related regularity techniques. The following theorem is a model result.
\begin{theorem}[\cite{dkk24}]\label{anit}
Let $u\in W^{1,1}_{\textnormal{loc}}(\Omega)$ be a minimizer of  
\begin{equation}\label{ani}
\scalebox{0.92}{  $ \displaystyle w\mapsto \int_{\Omega}\left(\snr{Dw}^{p}+\sum_{i=1}^{n}\snr{D_{i}w}^{q_{i}}\right)\dx,  \ \ 
    2\le p\le q_{1}\le \cdots\le q_{n}.$}
\end{equation}
Assume that
\begin{equation}\label{ani.1}
      q_{n}<\frac{p(n-1)}{n-3}\quad \mbox{if} \ \ n\ge 4
\end{equation}
and no other condition if $n=2,3$. 
Then $u\in W^{1,\infty}_{\textnormal{loc}}(\Omega)$.
\end{theorem}
\noindent Theorem \ref{anit} covers the models originally considered by Marcellini \cite{mar89,mar91}. If $p=2$ in \eqref{ani}, the bound in \eqref{ani.1} reduces to $q<2(n-1)/(n-3)$, which is precisely the threshold violated by \eqref{0.16}, so that Theorem \ref{anit} is sharp for polynomial-type integrals with quadratic growth from below. In two and three space dimensions no condition on $p,q$ is needed, as implicitly suggested by \eqref{ani.1}. Concerning the superlinear counterpart of  \eqref{0.13}, Choe \cite{cho92} and Esposito, Leonetti and Mingione \cite{elm99} showed different gradient regularity results for a priori locally bounded minimizers provided that $q<p+1$ and $q<p+2$, respectively, consistently with \eqref{0.13} when formally letting $\mu=2-p$.

\section{Uniformly elliptic Schauder estimates}\label{s3}
The focus of Schauder theory for elliptic equations or variational integrals is to quantify the effect of external data, i.e., coefficients, on the regularity of solutions. 
\subsection{Classical Schauder (a model case)}\label{s3.1}\vspace{1.5mm}
By Weyl's lemma, $L^1$-regular distributional solutions to the Laplace equation
$-\Delta u=0
$
are smooth. This easily extends to linear elliptic equations with constant coefficients.  The subsequent question is how much of this regularity is preserved when plugging in non-constant coefficients, or, more precisely, how the regularity of coefficients affects that of solutions.  Specifically, with $\tx{A}\colon \Omega\to\mathbb{R}^{n\times n}$ being a bounded and elliptic matrix, i.e., $\mathbb{I}_{n\times n}\approx \tx{A}$ in the sense of matrices, what can be said on the regularity of weak solutions\footnote{Although this is not strictly necessary in the linear case when coefficients are H\"older continuous, here we assume to deal with energy solutions, that is, distributional solutions that belong to the reference energy space $W^{1,2}(\Omega)$. These are usually called weak solutions.} to
\begin{equation}\label{0.17}
-\textnormal{div}\left(\tx{A}(x)Du\right)=0\qquad \mbox{in} \ \ \Omega\mbox{\ ?}
\end{equation}
Since $\tx{A}$ and $Du$ stick together in \eqref{0.17}, a natural guess is 
\begin{equation}\label{0.18}
\tx{A}\in C^{0,\beta}_{\textnormal{loc}}(\Omega,\mathbb{R}^{n\times n}) \ \Longrightarrow \ Du\in C^{0,\beta}_{\textnormal{loc}}(\Omega,\mathbb{R}^{n}),\vspace{0.2mm}
\end{equation}
which is in fact true for all $\beta\in (0,1)$. Results in the spirit of \eqref{0.18} were obtained by Hopf, Caccioppoli, Giraud and Schauder (1929-1934), including global versions. Later on, streamlined and different approaches were found by several other authors. All the methods available unavoidably exploit the quantitative information on the power-type decay of the modulus of continuity of coefficients (H\"older continuity), to show that energy solutions to \eqref{0.17} are close at all scales to harmonic-type maps, such as for instance their $\tx{A}(x_{0})$-harmonic lifting $v$ in $B_{r}$
\begin{equation*}
  \displaystyle
  \ -\textnormal{div}(\tx{A}(x_{0})Dv)=0\  \mbox{in} \  B_{r}, \quad  v=u \ \mbox{on} \  \partial B_{r}.
\end{equation*}
Indeed, ellipticity yields the homogeneous comparison estimate
\begin{equation}\label{0.19}
\mint_{B_{r}}\snr{Du-Dv}^{2}\dx\lesssim r^{2\beta}\mint_{B_{r} }\snr{Du}^{2}\dx.
\end{equation}
On the other hand, standard theory for linear elliptic equations with constant coefficients grants homogeneous decay estimates as 
\begin{equation}\label{0.20}
\scalebox{0.93}{  $ \displaystyle
\mint_{B_{\sigma}}\snr{Dv-(Dv)_{B_{\sigma}}}^{2}\dx\lesssim \left(\frac{\sigma}{\varrho}\right)^{2}\mint_{B_{\varrho}}\snr{Dv-(Dv)_{B_{\varrho}}}^{2}\dx,$}
\end{equation}
for all concentric balls $B_{\sigma}\subset B_{\varrho}\subset B_r$. Estimates \eqref{0.19}-\eqref{0.20} can be matched and iterated to deliver 
$$
\mint_{B_{r}}\snr{Du-(Du)_{B_{r}}}^{2}\dx \lesssim r^{2\beta}
$$
for all balls $B_r\Subset \Omega$, which implies the local $\beta$-H\"older continuity by certain integral characterization of H\"older continuity due to Campanato and Meyers.\footnote{The one described here is in fact Campanato's classical approach to Schauder estimates \cite{camp}.} This line of proof extends to $W^{1,p}$-regular distributional solutions to nonautonomous, quasilinear operators, such as
\begin{equation}\label{0.21}
-\textnormal{div}(\gamma(x)\snr{Du}^{p-2}Du )=0\qquad \mbox{in} \ \ \Omega,
\end{equation}
with $1\lesssim \gamma (\cdot)\in C^{0,\beta}_{\textnormal{loc}}(\Omega)$, $\beta \in (0,1)$ and $1<p<\infty$. This is due to the work of Manfredi \cite{man88}, Giaquinta and Giusti \cite{gg82,gg83,gg84}, and DiBenedetto \cite{dib83}. Also in this case $Du$ is locally H\"older continuous.\footnote{Gradient H\"older continuity of energy solutions holds but, in general, not with the sharp exponent $Du \in C^{0,\beta}$, due to the fact that the equation is degenerate.} Note that Schauder estimates obviously imply Lipschitz estimates 
\eqn{classicalsc}
$$
\mbox{$C^{1,\beta}$-estimates $\Longrightarrow C^{0,1}$-estimates}
$$
and this is in general the only way to get Lipschitz estimates when in presence of H\"older continuous coefficients as the equations considered cannot be differentiated. The key point the above techniques rely on is that all the a priori estimates involved, such as \eqref{0.19}-\eqref{0.20}, are homogeneous, and, as such, can be iterated. In turn, this is a feature of uniform ellipticity. When uniform ellipticity fails, a priori estimates are in general not homogeneous and these classical schemes fail as well. 
\subsection{More uniformly elliptic Schauder}\label{s3.2} \vspace{1.5mm}
The double phase functional\footnote{We denote $C^{\alpha}\equiv C^{[\alpha], \alpha -[\alpha]}$ when $\alpha$ is not an integer, and $[\alpha]$ denotes its integer part.}
\begin{equation}\label{0.22}
\begin{array}{c}
\displaystyle
w\mapsto \int_{\Omega}\left(\snr{Dw}^{p}+a(x)\snr{Dw}^{q}\right)\dx\vspace{1.5mm}\\ \displaystyle
1<p\le q,\quad 0\le a(\cdot)\in C^{\alpha}(\Omega)
\,, \quad   \alpha\in (0, \infty)
\end{array}
\end{equation}
was first considered by Zhikov in the context of homogenization of strongly anisotropic materials and of the study of Lavrentiev phenomenon \cite{jko94}. It only satisfies nonstandard growth conditions of $(p,q)$ type as in \eqref{soddi} but it is still uniformly elliptic in the sense that, with  $\tx{F}(x,z):= \snr{z}^p+a(x)\snr{z}^q$, the  ellipticity ratio $\mathcal{R}_{\tx{F}}(x,z)$ remains uniformly bounded. Indeed, Schauder-type results hold as in the following theorem.
\begin{theorem}[Baroni, Colombo and Mingione  \cite{bcm18, cm15}]\label{bcmt} 
Let $u\in W^{1,1}_{\textnormal{loc}}(\Omega)$ be a minimizer of \eqref{0.22}  with $\alpha \in (0,1]$. If\vspace{1mm}
\begin{itemize}
    \item either  $q/p\le 1+\alpha/n$, \vspace{0.5mm}
    \item or $u\in L^{\infty}_{\textnormal{loc}}(\Omega)$ and $q\le p+\alpha$,
    \end{itemize}
    \vspace{1mm}
    then $Du$ is locally H\"older continuous.\vspace{1.5mm}
\end{theorem}
\noindent The key of the proof is that the uniform ellipticity of the double phase functional \eqref{0.22} allows to implement a few more refined, nonstandard perturbation arguments. Specifically, recalling the discussion in Section \ref{s3.1}, minimizers to frozen functionals of the type
\begin{equation}\label{frozen}
w\mapsto \int_{B_{r}(x_{0})}\left(\snr{Dw}^{p}+a(x_{0})\snr{Dw}^{q}\right)\dx
\end{equation}
have locally H\"older continuous gradient and enjoy good reference estimates. This is in fact a consequence of the fact that functionals as in \eqref{frozen} are uniformly elliptic for every choice of $x_{0}$. On the other hand, the aforementioned nonstandard growth conditions, impacting solely on the comparison estimates, can be compensated via certain delicate schemes of reverse H\"older inequalities and higher integrability lemmas. Eventually, the approach of \cite{bcm18, cm15} was extended in \cite{bb25,ho22,ho22b} to treat large classes of uniformly elliptic integrals with nonstandard growth conditions as the double phase one; see also \cite{dm20} where the bound $q<p+\alpha$ is proved to be effective also in the vectorial case.

 \subsection{Soft nonuniform ellipticity and hard irregularity}\label{s3.22} \vspace{1.5mm}
 One might argue that Theorem \ref{bcmt} is incomplete as, being the double phase functional uniformly elliptic, Schauder-type results should hold with no restrictions on $p,q,\alpha$. Surprisingly enough, as first discovered in \cite{elm04}, the conditions imposed on such quantities in Theorem \ref{bcmt} are necessary and the result is sharp. 
 In fact, building on certain Zhikov's two-dimensional examples \cite{jko94} in the setting of the Lavrentiev phenomenon, in \cite{elm04, fmm04} a novel, sharp phenomenology was disclosed, revealing the failure of Schauder estimates in general, notwithstanding the uniform ellipticity of the problem considered. 
 \begin{theorem}[Fonseca, Mal\'y and Mingione \cite{fmm04}]\label{t4}
For every choice of parameters
\begin{equation}\label{0.25.1}
\begin{cases}
\displaystyle
\    1<p<n<n+\alpha<q<\infty,\quad  \alpha\in (0, \infty)
    \vspace{1mm}\\ \displaystyle
\  \qquad n\geq 2 , \qquad \varepsilon >0
 \end{cases}
\end{equation}
there exists a double phase integral \eqref{0.22}, 
a related minimizer $u\in W^{1,p}_{\loc}(\Omega)\cap L^{\infty}_{\loc}(\Omega)$, and a
closed set $\Sigma\subset\Omega$ with $\Dim(\Sigma) > n-p-\varepsilon,$ such that all the points of $\Sigma$ are non-Lebesgue points of
the precise representative of $u$.
\end{theorem}
\noindent In fact, in the same range \eqref{0.25.1}, non-$W^{1,q}$-regular, yet bounded minima of \eqref{0.22} with one-point singularity were constructed by Esposito, Leonetti and Mingione in \cite{elm04}. New constructions of singular minimizers eventually came, combining and improving the above features.
\begin{theorem}[Balci, Diening and Surnachev \cite{bds20,bds23}]\label{t42}
For every choice of the parameters 
\begin{equation}\label{dopo}
 q> p+\alpha \max \biggl\{1, \frac{p-1}{n-1}\biggr\}, \quad   \alpha\in (0, \infty), \quad p>1
\end{equation}
there exists a double phase integral \eqref{0.22}, and a related minimizer~$u\in W^{1,p}_{\loc}(\Omega)\cap L^{\infty}_{\loc}(\Omega)$, such that
$u\notin W^{1,d}$ for any~$d>p$. If $p<n$, there exists a closed set~$\Sigma\subset\Omega$ of non-Lebesgue points of $u$ with
$ \dim_{\mathcal{H}}(\Sigma)=n-p$.
\end{theorem}
\noindent The occurrence of irregular minima is not only explained in terms $(p,q)$-growth conditions. More is actually there, i.e., a softer form of nonuniform ellipticity hidden in \eqref{0.22}, that cannot be detected by using the classical ellipticity ratio $\mathcal{R}_{ \tx{F}}(x,z)$ in \eqref{0.8}, but rather considering a larger, nonlocal quantity accounting for the contribution of coefficients to the ellipticity of the functional over sets with positive measure.  Specifically, with $\BBB\subset \Omega$ being a ball, we consider the nonlocal ellipticity ratio \cite{dm21} defined as
\begin{equation}\label{0.24}
 \bar{\mathcal{R}}_{\tx{F}}(z,\BBB):=\frac{\sup_{x\in \BBB}\mbox{highest eigenvalue of} \ \partial_{zz}\tx{F}(x,z)}{\inf_{x\in \BBB}\mbox{lowest eigenvalue of} \ \partial_{zz}\tx{F}(x,z)}
\end{equation}
for $\snr{z}\not = 0$. Observe that $\mathcal{R}_{ \tx{F}}(x,z) \leq  \bar{\mathcal{R}}_{\tx{F}}(z,B)$ for $x\in \BBB$ and that the best upper bound obtainable on $ \bar{\mathcal{R}}_{\tx{F}}(z,B)$ is this time 
\begin{equation*}
 \bar{\mathcal{R}}_{\tx{F}}(z, \BBB)\lesssim_{p,q}1+\nr{a}_{L^{\infty}(B)}\snr{z}^{q-p}\,.
\end{equation*}
Moreover, if $a(\cdot)$ vanishes at some point in $\bar{\BBB}$, then 
$$\nr{a}_{L^{\infty}(B)} \snr{z}^{q-p} \lesssim  \bar{\mathcal{R}}_{\tx{F}}(z, \BBB) $$ so that $\bar{\mathcal{R}}_{\tx{F}}(z)\to \infty$ as $\snr{z}\to \infty$ if $a(\cdot)$ does not vanish identically in $\BBB$. This could be considered as a weaker form of nonuniform ellipticity, eventually generating singular minimizers although in presence of regular coefficients and classical uniform ellipticity.  Indeed, note that $ \bar{\mathcal{R}}_{\tx{F}}(z, \BBB)$ remains bounded when $a(\cdot)$ stays quantitatively away from zero on $\bar{\BBB}$, and in this case the same proof of Theorem \ref{bcmt} implies that $Du$ is locally H\"older continuous in $\BBB$, this time with no restriction on $p,q,\alpha$. 
\subsection{Fractal cones and malicious competitors}\vspace{1.5mm}

\begin{figure}[!ht]
  \centering
  \begin{tikzpicture}[scale=1.95]
    \node at (0,1.15) {Competitor~$u_{*}$};
    \draw[dashed] (-1,-1) -- (-1,+1) -- (+1,+1) -- (+1,-1) --cycle;
    \node at (0,.7) {\scalebox{0.8}{$m$}};
    \node at (0,-.7) {\scalebox{0.8}{$-m$}};
    \draw (-1,-1/2) -- (1,1/2);
    \draw (-1,1/2)-- (1,-1/2);
    \node at (+.6,0) {\scalebox{0.8}{$0$}};
    \node at (-.6,0) {\scalebox{0.8}{$0$}};
    \filldraw[pattern=vertical lines] (-1,-1/2) -- (0,0) -- (-1,-1/4);
    \filldraw[pattern=vertical lines] (-1,+1/2) -- (0,0) -- (-1,+1/4);
    \filldraw[pattern=vertical lines] (+1,-1/2) -- (0,0) -- (+1,-1/4);
    \filldraw[pattern=vertical lines] (+1,+1/2) -- (0,0) -- (+1,+1/4);
  \end{tikzpicture}
  \quad
  \begin{tikzpicture}[scale=1.95]
    \node at (0,1.15) {Coefficient $a(\cdot)$};
    \draw[dashed] (-1,-1) -- (-1,+1) -- (+1,+1) -- (+1,-1) --cycle;
    
    \node at (.6,0) {\scalebox{0.8}{$a=0$}};
    \node at (-.55,0) {\scalebox{0.8}{$a=0$}}; 
    \draw (-1/2,-1)--(+1/2,1);
    
    \draw(+1/2,-1)-- (-1/2,+1) ;
    \node at (0,.75) {\scalebox{0.8}{$a>0$}};
    \node at (0,-.75) {\scalebox{0.8}{$a>0$}};
    
    \filldraw[pattern=horizontal lines] (-1/2,-1) -- (0,0) -- (-1/4,-1);
    \filldraw[pattern=horizontal lines] (+1/2,-1) -- (0,0) -- (+1/4,-1);
    \filldraw[pattern=horizontal lines] (-1/2,+1) -- (0,0) -- (-1/4,+1);
    \filldraw[pattern =horizontal lines]
    (+1/2,+1) -- (0,0) -- (+1/4,+1);

  \end{tikzpicture}
  \caption{Competitor $u_{*}$ vs coefficient $a(\cdot)$. Figure \ref{fig} is a modification of the one in \cite{bds20}.}
  \label{fig}
\end{figure}
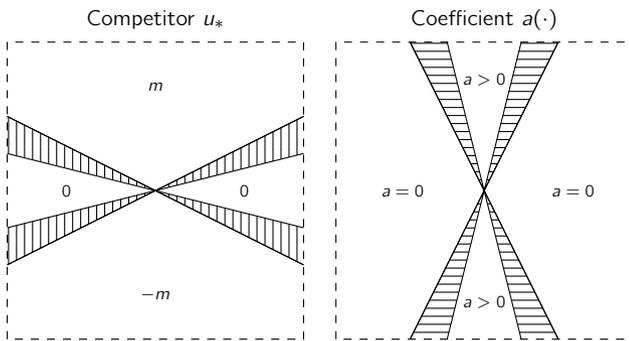

 The key to Theorems  \ref{t4}-\ref{t42} (we concentrate here on the second one, case $p<n$) is in the blending of three main ingredients.\vspace{1.5mm}
\begin{itemize}
    \item A merely $W^{1,p}$-regular map $u_{*}$ - the malicious competitor - attaining opposite values $m$ and $-m$ on the top and the bottom of $\Omega=[-1,1]^n$ and whose singularities can be distributed along a Cantor-type fractal $\mathcal{C}$ whose Hausdorff dimension $\dim_{\mathcal{H}}$ equals $n-p$. Here $m\ge1$ is a large constant.
    \item A Lipschitz-regular boundary datum $u_{0}$, with $u_0\equiv u_{*}$ on $\partial \Omega$.
    \item A nonnegative, $\alpha$-H\"older continuous coefficient $a(\cdot)$ vanishing where $\snr{Du_{*}}$ is positive, see Figure \ref{fig}.\vspace{1.5mm}
\end{itemize}
The last bullet point means that $a(x)\snr{Du_{*}}^{q}=0$ in $\Omega$, and therefore $u_{*}$ is a finite energy competitor in the Dirichlet problem driven by integral \eqref{0.22}, with $p,q,\alpha$ as in \eqref{dopo}, and boundary datum $u_{0}$. Basic Direct Methods of the calculus of variations yield the existence of a unique solution 
\begin{equation*}
u \mapsto \min_{w \in u_{0}+ W^{1,p}_{0}(\Omega)}  \int_{\Omega}\left(\snr{Dw}^{p}+a(x)\snr{Dw}^{q}\right)\dx
\end{equation*}
whose energy is set low, being controlled via minimality by the $p$-energy of $u_{*}$. Recalling $u_{0}-u_{*}\in W^{1,p}_{0}(\Omega)$ and that $u_{0}$ reaches opposite values on the top and the bottom of $\Omega$, a sufficiently large choice of $m$ ensures that the minimum $u$ "does not have enough energy" to cover the gap between the lower trace $-m$ and the upper one $m$ without developing discontinuities. In other, more accurate terms, a delicate combination of energy and trace estimates allows proving that $\Sigma_{u}$, the set of essential discontinuity points of $u$, contains a piece of fractal $\mathcal{C}$, thus forcing $\dim_{\mathcal{H}}(\Sigma_{u})=n-p$. This implies\footnote{The Hausdorff dimension of the set of non-Lebesgue points of a $W^{1,d}$-regular function does not exceed $n-d$, $d\leq n$.} that $u\not \in W^{1,d}_{\textnormal{loc}}(\Omega)$ for all $d>p$ showing that higher Sobolev regularity is in general unattainable under condition \eqref{dopo}. 
This construction is paradigmatic of the idea that, once identified the right (bad) competitor, and a related geometry of the coefficient, minimality can be used to produce singularities rather than proving regularity properties. The strength of these examples lies in the following aspects:\vspace{1.5mm}
 \begin{itemize}
 \item Minimizers, which are simply as bad as any other competitor. 
     \item Scalar setting. This is a genuinely nonstandard growth conditions phenomenon, in contrast with standard cases, where to produce singularities one needs to look at vectorial problems \cite{deg68} or to violate the initial energetic information \cite{serrin}. 
     \item No degeneracy issues \cite{ura68}. The integrand can be further made nondegenerate by replacing $\snr{Dw}$ with $(\snr{Dw}^2+1)^{1/2}$.  
     \item Lipschitz domains and boundary data. 
 \end{itemize}
\section{Schauder estimates and $\mu$-ellipticity}\label{finalsec}
The perturbation-based circle of ideas and techniques discussed in Section \ref{s3.2}  breaks down when genuine nonuniform ellipticity is involved: both reference and comparison estimates become nonhomogeneous, and the perturbative approaches, based on iterations, become unviable. The validity of Schauder estimates in the nonuniformly elliptic setting was a longstanding open problem raised at various stages in the literature: see e.g., \cite[Page 7]{iva84} on classical results of Ladyzhenskaya and Ural'tseva \cite{lu68}, Giaquinta and Giusti's paper \cite{gg84}, and  its MathSciNet review\footnote{Math. Rev. MR0749677.} by Lieberman. A complete solution eventually appeared in \cite{dm23a,dm25} in the $(p,q)$-setting, and in \cite{dm23b,ddp24} in the nearly linear, $\mu$-elliptic one. The novel techniques introduced in \cite{dm23a,dm23b,dm25,ddp24} reverse  the classical paradigm in \eqref{classicalsc} to obtain gradient estimates when dependency on coefficients is H\"older continuous. Indeed, for the first time gradient $L^\infty$-bounds are not derived as a consequence of $C^{1,\beta}$-bounds (in turn obtained via perturbation), but are rather derived directly, and eventually used to prove $C^{1,\beta}$-estimates. In other words, we return to \eqref{classicalsc0} although the functionals and the equations considered here are non-autonomous and non-differentiable. We shall try to give an overview of some of the ideas leading to establish Schauder estimates for certain classes of functionals with nearly linear growth. As explained immediately after Theorem \ref{t2} and displayed in \eqref{classicalsc0}, we can concentrate on Lipschitz estimates. 

\subsection{Nearly linear Schauder and intrinsic Bernstein functions} \vspace{1.5mm}The main models to initially  keep in mind are the logarithmic energies \eqref{0.6}-\eqref{0.7}, but now  also featuring H\"older continuous coefficients. 
\begin{theorem}[\cite{dm23b}]\label{t4.1}
Let $u\in W^{1,1}_{\textnormal{loc}}(\Omega)$ be a minimizer of functional 
\begin{equation}\label{0.26}
    w\mapsto \int_{\Omega}\left(\gamma(x)\snr{Dw}\log(1+\snr{Dw})\right)\dx
\end{equation} 
where $1\lesssim \gamma(\cdot)\in C^{0,\beta}(\Omega)$, $\beta\in (0,1)$. 
Then $u\in C^{1,\beta}_{\textnormal{loc}}(\Omega)$\footnote{In fact in \cite{dm23b} we proved that $u\in C^{1,\beta/2}_{\textnormal{loc}}(\Omega)$, but the improvement to the full exponent $u\in C^{1,\beta}_{\textnormal{loc}}(\Omega)$ can be easily reached, see \cite{ddp24}, by arguing as in \cite[Section 5]{dm25}. The main point in Theorem \ref{t4.1} is as usual to get that $Du\in L^{\infty}_{\loc}$, although the adaption of the standard perturbation methods to get gradient H\"older continuity once $Du$ is known to be locally bounded still requires care.}.
\end{theorem} 
\noindent Analogous results hold if in \eqref{0.26} the $L\log L$ integrand is replaced by the iterated logarithmic one in \eqref{0.6.1}. Theorem \ref{t4.1} is actually a special case of a more general result covering nonautonomous $\mu$-elliptic functionals, like for instance those exhibited by nearly linear double phase integrals of the type
\begin{equation}\label{0.27}
\scalebox{0.87}{
$\begin{array}{c}
\displaystyle
    w\mapsto \int_{\Omega}\left(\snr{Dw}\log(1+\snr{Dw})+a(x)(\snr{Dw}^{2}+s^{2})^{q/2}\right)\dx\vspace{1.5mm}\\
    \displaystyle
    1<q,\quad 0\le a(\cdot)\in C^{0,\alpha}(\Omega),\quad 0\leq s \leq 1.
    \end{array}$}
\end{equation}
This can be considered as the borderline configuration of \eqref{0.22} as $p\to 1$, while actually approaching nearly linear growth conditions. A key point here is that, on the contrary of \eqref{0.22}, the functional in \eqref{0.27} is not uniformly elliptic. This is easily seen at those points $x$ where $a(x)=0$, where the integrand reduces to $\snr{Dw}\log(1+\snr{Dw})$, which is nonuniformly elliptic.  Nevertheless, also in this case it is possible to achieve maximal regularity for minima, and under optimal structural conditions regulating nonuniform ellipticity. Furthermore this extends to larger classes of functionals, to which \eqref{0.27} belongs to, of the type 
\begin{equation*}
 w\mapsto \int_{\Omega} \left(\gamma(x)\mathfrak{L}(Dw)+a(x)(\snr{Dw}^{2}+s^{2})^{q/2}\right)\dx
\end{equation*}
where $\gamma (\cdot)$ is as in Theorem \ref{t4.1} and 
\begin{equation*}
\frac{\snr{\xi}^{2}}{(\snr{z}^{2}+1)^{\mu/2}}\lesssim \langle\partial_{zz}\mathfrak{L}(z)\xi,\xi\rangle\lesssim \frac{\snr{\xi}^{2}}{(\snr{z}^{2}+1)^{1/2}}
\end{equation*}
is assumed to hold for $\mu\approx 1$ (i.e., for $1\leq \mu< \mu_m\equiv \mu_m(n,q,\alpha)<2$) and every choice of $z,\xi \in \mathbb R^n$, see \cite{dm23b} for details. For instance, all models featuring iterated logarithms as 
\begin{equation}\label{0.28}
    w\mapsto \int_{\Omega}\left(\tx{L}_{i+1}(Dw)+a(x)(\snr{Dw}^{2}+s^{2})^{q/2}\right)\dx
\end{equation}
are included, where $s\in [0,1]$, and $a(\cdot)$ and $q$ as in \eqref{0.27}. In this respect, the following holds.
\begin{theorem}[\cite{dm23b,ddp24}]\label{t4.2}
Let $u\in W^{1,1}_{\textnormal{loc}}(\Omega)$ be a minimizer of functionals in \eqref{0.27} or in \eqref{0.28} with $\alpha \in (0,1)$. If 
\begin{itemize}
    \item either $q<1+\alpha/n$,\vspace{0.5mm}
    \item or $u\in L^{\infty}_{\textnormal{loc}}(\Omega)$ and 
        $q<1+\alpha$,
\end{itemize}
    then $Du$ is locally H\"older continuous in $\Omega$.
Moreover, in the nonsingular case $s>0$, we have $u\in C^{1,\alpha}_{\textnormal{loc}}(\Omega)$ for \eqref{0.27}, and $u\in C^{1,\alpha/2}_{\textnormal{loc}}(\Omega)$ for \eqref{0.28}.
\end{theorem}
\noindent The bound $q<1+\alpha$ reveals to be sharp, as the following holds true.
\begin{theorem}[\cite{ddp24}]\label{t4.3}
For every choice of the parameters $\alpha,\varepsilon> 0$, $q>1$, such that $q>1+\alpha$, $0<\varepsilon<\min\{q-1-\alpha, n-1\}$ there exists a double phase integral \eqref{0.27} and related minimizer $u \in W^{1,1}_{\loc}(\Omega)\cap L^{\infty}_{\loc}(\Omega)$ such that $u\not \in W^{1,p}_{\loc}(\Omega)$ for all $p>1+\varepsilon$. In particular, the Hausdorff dimension of the set of non-Lebesgue points of the precise representative of $u$ is at least equal to $n-1-\varepsilon$.
\end{theorem}
\noindent Although the outcome is formally the same, the approach to Theorem \ref{t4.2} is completely different from the one of Theorem \ref{bcmt}. Indeed, while the functional in \eqref{0.22} is {\em uniformly elliptic}, the functionals in \eqref{0.27} and \eqref{0.28} are not. Therefore, perturbation approaches
of the type considered in \cite{bb25,bcm18,cm15,ho22,ho22b} fail to deliver results and more complicated, completely different routes are necessary.  In this respect, note that the functionals in \eqref{0.27}-\eqref{0.28} globally satisfy assumptions \eqref{0.5} for any $\mu>1$ (actually we can take $\mu=1$ in case \eqref{0.27}), but the simple appeal to such properties is not sufficient to prove Theorem \ref{t4.2} and more is needed. There are in fact three main points in the proof of Theorem \ref{t4.2} and in the following we shall restrict for simplicity to the case where the bound $q<1+\alpha$ is considered.  A first key idea is to fully exploit the specific structure of the functional to rebalance the significant loss of ellipticity due to degenerate nonuniform ellipticity. This is achieved via a novel, intrinsic version of the Bernstein technique combining fractional estimates and nonlinear potential theoretic methods. In the uniformly elliptic case of functionals of the type \eqref{plap} one observes that a function of the type 
\begin{equation}\label{intri}
v(x)=(\snr{Du(x)}^2+1)^{p/2}
\end{equation} is a subsolution to a linear, uniformly elliptic equation, and, as such, it is bounded. This follows from the possibility of differentiating the related Euler-Lagrange equation and the fact that the functional is uniformly elliptic. Both things fail for \eqref{0.28}. Indeed, recall that Euler-Lagrange equation to \eqref{0.28} is 
\begin{equation*}
\scalebox{0.91}{ $
-\textnormal{div} \, (\partial_z\tx{L}_{i+1}(Du))  -q\textnormal{div}
(a(x)( s^{2}+\snr{Du}^{2})^{(q-2)/2}Du)=0$}
\end{equation*}
and therefore is it not differentiable by the H\"older continuity of the coefficient $a$. The idea is then to replace the function $v$ in \eqref{intri} by another, more intrinsic Bernstein function, incorporating larger information on the structure of the integrand and its ellipticity, namely
\begin{eqnarray*}
  \notag  \tx{E}(x) &:= & \frac{1}{2-\mu}\Big[(\snr{Du(x)}^{2}+1)^{1-\mu/2}-1\Big]\\[8pt]
    && \quad +(1-1/q) \tx{a}(x)\Big[(\snr{Du(x)}^2+s^2)^{q/2}-s^{q}\Big]
\end{eqnarray*}
where in fact $\mu \in [1,2)$ is the one for which  \eqref{0.5} is satisfied. In turn, this function is shown to satisfy a renormalized, fractional Caccioppoli-type inequality,\footnote{Of course \eqref{0.30} makes sense as an a priori estimate and must be fixed via an approximation argument where original minimizers are the limit of minima of certain more regular, uniformly elliptic functionals.} i.e., 
\begin{eqnarray}\label{0.30}
\notag  r^{2\beta}[(\tx{E}-\kappa)_{+}]_{\beta,2;B_{r/2}}^{2} & \lesssim  & \texttt{M}^{2\tx{b}_{1}}\int_{B_{r}}(\tx{E}-\kappa)_{+}^{2}\dx\vspace{2.5mm}\\[8pt]
&&  \hspace{-9mm}+\texttt{M}^{2\tx{b}_{2}}r^{2\alpha}\int_{B_{r}}(\snr{Du}^{m}+1)\dx
\end{eqnarray}
holds 
for any $\kappa\ge 0$, all balls $B_{r}\subset \Omega$ with radius $r$, suitable numbers $\beta\in (0,\alpha)$, $\tx{b}_{1},\tx{b}_{2},m\geq 1$ and $\tx{M}$ such that $\tx{M}\geq \nr{Du}_{L^{\infty}(B_{r})}$. On the left-hand side in \eqref{0.30} there appears the classical fractional Gagliardo norm, which is defined as 
\begin{equation*}
[v]_{\beta,2;A}^{2}:= \int_A\int_A \frac{\snr{v(x)-v(y)}^2}{\snr{x-y}^{n+2\beta}}\dx \dy 
\end{equation*}
whenever $A\subset \mathbb{R}^{n}$ is an open set and $v\colon A \to \mathbb R$ is a measurable function. The term renormalized accounts for the fact that inequality \eqref{0.30} is homogeneous with respect to $\tx{E}$, despite integrals \eqref{0.27}-\eqref{0.28} are not. This is exactly the feature allowing to apply the nonlinear potential machinery mentioned above. The price to pay is the appearance of multiplicative constants depending on $\nr{Du}_{L^{\infty}(B_{r})}$ (via $\texttt{M}$). Such constants must be carefully kept under control all over the proof and reabsorbed at the very end. This will be a point where the bound $q< 1+\alpha$ assumed in Theorem \ref{t4.2} is used in a crucial manner. The validity of \eqref{0.30} is established via a nonlinear dyadic/atomic decomposition technique, finding its roots in \cite{km05}, that resembles the one used for Besov spaces in the setting of Littlewood-Paley theory. Fractional Caccioppoli inequalities of the type in \eqref{0.30}, first pioneered in \cite{Mipo} in the setting of Nonlinear Potential Theory, eventually allow to prove boundedness of $\tx{E}$ via a nonlinear potential theoretic version of De Giorgi's iteration, that made its first appearance in \cite{km93}. In this respect, here a more delicate and quantitative form of such ideas is needed \cite{dm23a}. The boundedness of $\tx{E}$ obviously implies the one of $Du$. Back to the proof of the boundedness of $\tx{E}$, we point out that the nonlinear potentials used in the estimates are of the type introduced by Khavin and Maz'ya \cite{hm72} and deeply studied by Adams, Hedberg, Meyers, Wolff. Specifically, these are of the form
\begin{equation*}
\mathbf{P}_{\sigma}^{\vartheta}(f;x,r):=\int_{0}^{r}\varrho^{\sigma}\Big(\frac{1}{\snr{B_{\varrho}(x)}}\int_{B_{\varrho}(x)}\snr{f}\dy\Big)^{\vartheta}\frac{\textnormal{d}\varrho}{\varrho}
\end{equation*}
for parameters $\sigma,\vartheta>0$ and $f$ being an $L^1(B_{r}(x))$-regular vector field. Suitable choices of $\sigma$ and $\vartheta$ give back the standard Riesz potential $\mathbf{I}_{1}$ and the Wolff potential $\mathbf{W}_{1,p}$ \cite{km18}. The mapping properties of potentials among function spaces are known. Specifically,
\begin{equation}\label{0.32pp}
\nr{\mathbf{P}_{\sigma}^{\vartheta}(f;\cdot,r)}_{L^{\infty}(B_{\varrho})}\lesssim \nr{f}_{L^{m}(B_{r+\varrho})}^{\vartheta},
\end{equation}
holds whenever $n\vartheta>\sigma$, $m>n\vartheta/\sigma$ and $B_{r+\varrho}\subset \Omega$, \cite{dm23a, dm25}. In contrast with previous foundational contributions \cite{km93,km18}, where potentials are employed as ghosts of the representation formula to derive optimal regularity of solutions from data, in \cite{dm23a,dm23b,dm25,ddp24} potentials fit the fractional nature of \eqref{0.30} and sharply quantify how the rate of Hölder continuity of coefficients interacts with the growth of the terms they stick to towards $L^\infty$-estimates. Another main idea in this setting is a fractional Moser's iteration in Besov spaces already employed for problems with polynomial growth in \cite{dm25},\footnote{See \cite{bls18,dom04} for similar Besov spaces techniques in the context of degenerate integro-differential equations.} and which allows reading the H\"older continuity of the coefficient $a(\cdot)$ as fractional differentiability. This allows to gain arbitrarily high gradient integrability and therefore, in a sense, to quantitatively reduce the rate of nonuniform ellipticity of \eqref{0.27}-\eqref{0.28}. However, further  obstructions arise due to the severe loss of ellipticity in integrals at nearly linear growth. These require a limiting version of the aforementioned fractional Moser's iteration in \cite{dm25} yielding hybrid reverse H\"older inequalities of the form
\begin{equation}\label{0.31}
\scalebox{0.94}{  $
\nr{Du}_{L^{t}(B_{r/2})}\lesssim \texttt{M}^{\omega}\left(1+\frac{\nr{u}_{L^{\infty}(B_{r})}}{r}\right)^{\texttt{b}}\left(1+\nr{Du}_{L^{1}(B_{r})}^{1/t}\right)$}
\end{equation}
which is valid for all $1\le t<\infty$, $\omega\in (0,1)$, some $\texttt{b}>0$, any ball $B_{r}\subset \Omega$, and, as in \eqref{0.30}, and again $\tx{M}\geq \nr{Du}_{L^{\infty}(B_{r})}$. Estimate \eqref{0.31} is a sort of borderline interpolation inequality, and, when combined with \eqref{0.30}, allows working under the maximal ellipticity range $q<1+\alpha$. Once again, the price to pay is the appearance in the bounding constants of $\texttt{M}^{\omega}$, with $\omega$ that can be picked to be arbitrarily small,  to compensate the loss of ellipticity and trade between an arbitrarily high power of the modulus of the gradient and its $L^{1}$-norm.

\begin{remark}[Obstacles]
The techniques devised for Theorem \ref{t4.2} are flexible enough to deal with variational obstacle problems. Specifically, they allow to bypass the classical linearization procedure pioneered by Duzaar and Fuchs \cite{duz87,df86} and to prove gradient regularity in nondifferentiable, nonuniformly elliptic variational inequalities. Duzaar and Fuchs's approach turns constrained minimizers of homogeneous integrals into unconstrained minima of forced functionals whose right-hand side is a function of the second derivatives of the obstacle and of the gradient of coefficients. This again requires that $x\mapsto \partial_z\tx{F}(x,\cdot)$ is differentiable, which is not the case in the present setting. Alternative techniques, as those used by for instance Choe \cite{cho91}, only work in the uniformly elliptic case. On the other hand, the scheme supporting Theorem \ref{t4.2}, based on fractional differentiation and use of nonlinear potentials, can be tailored to account for obstacles in order to deliver sharp results also in the constrained case. For this we refer to \cite{ddp24}.
\end{remark}

\subsection{Back to the beginning}\label{s4.1}\vspace{1.5mm} We finally highlight a few formal connections between results of the type in Theorems \ref{t4.2}-\ref{t4.3} and some classical counterexamples to regularity for linear growth functionals 
constructed by Giaquinta, Modica and Souček \cite{gms79}. Consider the generalized\footnote{By "generalized" we mean that problem \eqref{gen1} must be extended to $BV$ and the functional appearing in \eqref{gen1} is actually replaced by a suitable relaxed form in which boundary data are penalized, \cite{bs13,gms79}.} Dirichlet problem involving the area functional in one dimension 
\begin{equation}\label{gen1}
\begin{cases}
\displaystyle
\ w\mapsto \min_{}\int_{-1}^{1}\sqrt{1+\gamma(x)\snr{ w'}^{2}}\dx\vspace{1mm}\\ \displaystyle
\ w(-1)=-w_{0},\qquad \quad w(1)=w_{0},
 \end{cases}
\end{equation}
where $\gamma(x):=1+x^{2}\left(\log(2/\snr{x})\right)^{4}$ and $w_{0}>0$ satisfies
\begin{equation*}
\infty > w_{0}>\int_{-1}^{1}\frac{1}{\sqrt{\gamma(x)-1}} \dx. 
\end{equation*}
\noindent The minimizer has a jump at zero, making $W^{1,1}$-regularity fail. The function $\gamma$ is not $C^{2}$-regular. 
In contrast, $C^{2}$-coefficients guarantee the possibility of a priori estimates \cite{lu70} in the style of \eqref{0.2}. 
This situation resembles the one of Theorems \ref{t4}, \ref{t42} and \ref{t4.3}. In this respect, the construction of Theorem \ref{t4.3} extends to linear growth double phase integrals such as
\begin{equation}\label{0.34}
w\mapsto \int_{\Omega}\left((1+\snr{Dw}^{m})^{1/m}+a(x)\snr{Dw}^{q}\right)\dx
\end{equation}
with $1<m,q$ and $0\leq a(\cdot)\in C^{\alpha}(\Omega)$.\footnote{Again, as above, one has to interpret this in a suitably relaxed way, considering competitors in $BV$ and a relaxed from of the functional.} This means that there is no hope of pointwise gradient regularity for minima of \eqref{0.34} whenever $q>1+\alpha$. An important point is that this  information comes only from the growth conditions of the integrand, and not from its ellipticity, i.e., the growth of the eigenvalues of the second derivatives. On the other hand, the bounds relevant in order to prove a priori estimates come from conditions on second derivatives like \eqref{0.5} or \eqref{0.14}. While these scale accordingly to the growth conditions of the integrand in superlinear growth regimes - like in the $(p,q)$ case, compare \eqref{0.14} and \eqref{soddi} - there might be a detachment when approaching the linear case: the integrand keeps growing linearly, but its derivatives can decay very fast (consider \eqref{0.6.2} with large $m$). In view of Theorem \ref{t4.2}, and proceeding formally, the condition for a priori regularity gradient estimates looks as $q < 2-\mu +\alpha$. In the case of anisotropic area-type functionals as in \eqref{0.34}, it is $\mu=m+1$. Coupling this with $q>1$ leads to  $\alpha >q+m-1$, which, for $m$ close to two, matches the need of $C^2$-regular coefficients in classical papers as \cite{lu70} and also aligns with counterexample \eqref{gen1}.  On the other hand, further imposing the restriction $\alpha\in (0,1)$ and then recalling that it must be $q<1+\alpha$, leads to $\mu\approx 1$, exactly as considered in \cite{dm23b,ddp24}, cf. Theorem \ref{t4.2}. This suggests that functionals \eqref{0.27}-\eqref{0.28} might be the limiting configurations for the validity of Schauder theory in presence of H\"older coefficients, convex anisotropy and $\mu$-ellipticity.

\begin{ack}
The author thanks Professor Anna Balci for sharing the drawings in Figure \ref{fig} and Professor Lorenzo Brasco for comments on a preliminary version. This work is supported by the European Research Council, through the ERC StG project NEW, nr.~101220121, and by the University of Parma through the action "Bando di Ateneo 2024 per la ricerca".
\end{ack}
\begin{authorinfo}
Cristiana De Filippis is full professor at the University of Parma, Italy. Her research focuses on the (ir)regularity of minima of possibly nonconvex multiple integrals, solutions to nonlinear elliptic or parabolic PDEs and integro-differential equations. Outside mathematics, she enjoys the company of her horse Tuono and her cat Anakin. 
\email{cristiana.defilippis@unipr.it}

\end{authorinfo}

\end{document}